\documentclass{article}

\usepackage{amsmath}
\usepackage{graphicx}
\newtheorem{theorem}{Theorem}
\newtheorem{lemma}{Lemma}

\title{Regularity Conditions for Convergence of Linear Statistics of GUE}
\author{Phil Kopel}
\date{}

\begin{document}
\maketitle

\begin{abstract}
We establish a central limit theorem for the unnormalized linear statistic of the Gaussian Unitary Ensemble under optimal conditions: the linear statistics converges if and only if the expression for the limiting variance is finite.
\end{abstract}

\section{Introduction}

Linear statistics are an indispensible tool for understanding the behavior of the spectra of large random matrices, as they describe the fluctuations of eigenvalues about their expected configurations. Indeed, these statistics can be thought of as  random matrix theory analogues of the central limit theorems of classical probability theory. 

The set-up is as follows: fix some test function $f$ and matrix $M_n$ with eigenvalues $\lambda_1,...,\lambda_n$. The associated linear statistic $N_n[f]$ then is defined to be:

\begin{eqnarray}
N_n[f]=\sum_{j=1}^n f(\lambda_j)-\mathbf{E}\left[\sum_{j=1}^n f(\lambda_j)\right]
\end{eqnarray}

The distribution of eigenvalues is frequently highly regular: for many random matrix models the linear statistic does not need to be normalized to converge to a limiting distribution, unlike sums of independent random variables, which experience fluctuations on the order of $\sqrt{n}$.

Much is known about these statistics already, however for many important matrix models the available results do not yet appear to be optimal (in the sense that they make unnecessarily strict demands on the regularity of test functions). Here, we will consider the case of Wigner matrices, which are a central object of study in random matrix theory and enjoy numerous applictions to both pure and applied mathematics \cite{Meh}.

We begin with some elementary definitions. A Wigner matrix $W_n$ is an $n$ dimensional random Hermitian matrix such that each pair of complex conjugate off-diagonal entries ($w_{ij}$ and $w_{ji}$, for $i\neq j$) and each diagnonal entry $w_{ii}$ are given by independent complex random variables with mean zero and variance $\frac{1}{n}$. Additionally, we will assume that the entries share third and fourth moments, and also obey the following estimate for some $c>0$:

\begin{eqnarray*}
\mathbf{P}(|w_{i,j}|\geq t^c)\leq e^{-t}
\end{eqnarray*}

The former assumption is explained by the sensitivity of the linear statistic to these moments, the latter assumption is a common technical restriction and is usually refered to as Condition C0 in the literature.

By the famous Wigner semicircular law, the eigenvalues $\lambda_1,...,\lambda_n$ of $W_n$  converge in probability to the semicircular distribution:
\begin{eqnarray*}
\lim_{n\to\infty}\mathbf{P}\left[\left|\frac{1}{n}\sum_{j=1}^n f(\lambda_j) - \frac{1}{2\pi}\int_{-2}^{2} f(x)\sqrt{4-x^2}dx\right|>\epsilon\right]=0
\end{eqnarray*}

For sufficiently nice test functions $f$, the limiting distribution of the associated linear statistic is that of a normal random variable with the following variance:
\begin{eqnarray*}
V_{\emph{WIG}}[f]=\frac{1}{4\pi^2}\int_{-2}^2\int_{-2}^2 \left(\frac{f(x)-f(y)}{x-y}\right)^2\frac{4-xy}{\sqrt{4-x^2}\sqrt{4-y^2}}dxdy\\
+\frac{\kappa_4}{4\pi^2}\left(\int_{-2}^2 f(x)\frac{2-x^2}{\sqrt{4-x^2}}\right)^2
-\frac{\kappa_2-2}{4\pi^2}\left(\int_{-2}^2 f(x)x\frac{2-x^2}{\sqrt{4-x^2}}\right)^2
\end{eqnarray*}

Here, $\kappa_2$ and $\kappa_4$ are the second and fourth cumulants of the entries, respectively. That some regularity ought to be needed from $f$ is already obvious: after all, the limiting variance clearly diverges if $f$ is discontinuous. Notice that continuity alone is not enough to prevent the limiting variance from diverging, however.

The exact regularity restrictions on $f$ have been improved by a long series of papers, a representitive sample of which includes Lytova and Pastur \cite{LP} (which uses Stein's method and requires 5 classical derivatives),  Scherbina \cite{Sch} (which uses a martingale approach and requires membership in the Sobolev class $H^{3/2+\epsilon}$), and most recently Sosoe and Wong \cite{SW} (which uses Littlewood-Paley theory and requires membership in the Sobolev class $H^{1+\epsilon}$).

If all entry distributions associated with a Wigner matrix are given by Gaussian random variables, the matrix is said to belong to the Gaussian Unitary Ensemble (GUE), a very important subclass of Wigner matrices and of random matrices in general. The advantage of working with these matrices is that their statistics can be computed explicity using a variety of classical formulas, the disadvantage of this advantage is that the methods employed are therefore extremely unlikely to be at all useful for other models.

The linear statistics of GUE are well-studied, with the strongest result currently in the literature, also due to Sosoe and Wong, ensuring a Gaussian limiting distribution for test functions with H\"older constant better than one half \cite{SW}. However, it is reasonable to suspect that the optimal condition ought to be the finiteness of the variance of the limiting distribution. In fact, this turns out to be exactly the case:

\begin{theorem}
\label{clt:gue}
Let $M_n$ be a GUE random matrix, and $f$ any bounded real test function. The linear statistic $S_n(f)$ converges in distribution to a random variable with finite variance if any only if the following expression is finite:
\begin{eqnarray}
V_{\emph{GUE}}[f]=\frac{1}{4\pi^2}\int_{-2}^2\int_{-2}^2 \left(\frac{f(x)-f(y)}{x-y}\right)^2
\frac{4-xy}{\sqrt{4-x^2}\sqrt{4-y^2}}dxdy
\end{eqnarray}

The limiting distribution is that of a normal random variable with zero mean and the above variance.
\end{theorem}

Notice that the variance expression is a simplified version of the variance in the Wigner case (since, for GUE,  all of the higher cumulants of the entry distributions vanish identically).  The proof is based on the explicit formulas available for the correlations functions of GUE. The key feature of our approach is that we will work with the expression for the variance directly, and avoid introducing various function space norms altogether.

\section{Background on GUE}

Here, we gather some standard background about GUE matrices which we will require. The reference for all of this material is the classical text \cite{AGZ}.

First, define the $n$-th (probabilist's) Hermite polynomial as follows:

\begin{eqnarray*}
h_n(x)=(-1)^ne^{x^2/2}\frac{d^n}{dx^n}e^{-x^2/2}
\end{eqnarray*}

These are polynomials with leading coefficient one which satisfy the following orthogonality relationship (with $\delta$
being the Kroenecker delta):

\begin{eqnarray*}
\int_{-\infty}^{\infty} h_{m}(x)h_{n}(x)e^{-x^2/2}dx=\sqrt{2\pi}n!\delta_{n,m}
\end{eqnarray*}

Define also the $n$-th normalized oscillator wave function:
\begin{eqnarray*}
\psi_n(x)=\frac{e^{-x^2/4}h_n(x)}{\sqrt{n!}(2\pi)^{1/4}}
\end{eqnarray*}

These functions satisfy the following identities:

\begin{eqnarray*}
x\psi_n(x)=\sqrt{n+1}\psi_{n+1}(x)+\sqrt{n}\psi_{n-1}(x)\\
\psi^{'}_n(x)=-\frac{x}{2}\psi_n(x)+\sqrt{n}\psi_{n-1}(x)\\
\psi^{''}_n(x)=-(n+\frac{1}{2}-\frac{x^2}{4})\psi_n(x)\\
\int_{\mathcal{R}} \psi_k(x)\psi_{j}(x)=\delta_{k,j}
\end{eqnarray*}

Further, define the kernel:
\begin{eqnarray*}
K_n(x,y)=K_n(y,x)=\sum^{n-1}_{j=0}\psi_j(x)\psi_j(y)
\end{eqnarray*}

$K_n(x,y)$ satisfies the following:
\begin{eqnarray*}
\int_{\mathcal{R}}K_n(x,x)dx=n\\
K_n(x,y)=\int_{\mathcal{R}}K_n(x,z)K_n(z,y)dz\\
K_n(x,y)=\sqrt{n}\frac{\psi_n(x)\psi_{n-1}(y)-\psi_{n-1}(x)\psi_n(y)}{(x-y)}
\end{eqnarray*}

This last identity  is the well known Christoffel-Darboux formula. The point of all this development is that $k$-point correlation functions, $p_k$, of the GUE ensemble are expressible in terms of these quantities by the Gaudin-Mehta formula:
\begin{eqnarray*}
p_k(\lambda_1,...,\lambda_k)=\det(K_n(\lambda_i,\lambda_j))_{1\leq i,j\leq n}
\end{eqnarray*}

Using the Christoffel-Darboux formula and the explicit formulas for one and two point correlation functions, it is straightforward to compute the following formula for the variance of the linear statistic:
\begin{eqnarray*}
\mbox{Var}(N_n[f])=\frac{n}{2}\int_{-2}^2\int_{-2}^2 (f(x)-f(y))^2
K_n(x,y)^2 dxdy\\
=\frac{n}{2}\int_{-2}^2\int_{-2}^2 \left(\frac{f(x)-f(y)}{x-y}\right)^2
\Psi_n(x,y) dxdy
\end{eqnarray*}

Here, we have defined:
\begin{eqnarray*}
\Psi_n(x,y)=\psi_n(\sqrt{n}x)^2\psi_n(\sqrt{n}y)^2-
\psi_n(\sqrt{n}x)\psi_{n-1}(\sqrt{n}x)\psi_n(\sqrt{n}y)\psi_{n-1}(\sqrt{n}y)
\end{eqnarray*}

In addition this classical material, we will also make use of the following bound:

\begin{lemma}
\label{boundbulk}
Fix $\delta>0$. There exists a constant $C$ such that the following estimate holds in the interval $[-2+\delta,2-\delta]$ (uniformly in $n$):
\begin{eqnarray*}
\sqrt{n}\max\left[\psi_n(\sqrt{n}x)^2,|\psi_{n-1}(\sqrt{n}x)\psi_n(\sqrt{n}x)|\right]\leq C
\end{eqnarray*}
The constant $C$ is allowed to depend on $\delta$.
\end{lemma}

This can easily be proven using, for instance, the machinery of Plancharel-Rotach asymptotics. An alternate and more elementary proof, which relies instead on a Tracy-Widom type estimate, is supplied in the Appendix for the interested reader.

Another lemma which we will need, whose proof is also in the Appendix, is:

\begin{lemma}
\label{conv}
Suppose $f\in L^{\infty}(\mathcal{R})$. Then:
\begin{eqnarray*}
\lim_{n\to\infty} \sqrt{n}\int f(x)\psi_n(\sqrt{n}x)^2dx=\frac{1}{\pi}\int \frac{f(x)}{\sqrt{4-x^2}}dx\\
\lim_{n\to\infty} \sqrt{n}\int f(x)\psi_n(\sqrt{n}x)\psi_{n-1}(\sqrt{n}x)dx=\frac{1}{\pi}\int \frac{xf(x)}{2\sqrt{4-x^2}}dx
\end{eqnarray*}
Furthermore, suppose $H(x,y)\in L^{\infty}[\mathcal{R}\times \mathcal{R}]$. Then:
\begin{eqnarray*}
\lim_{n\to\infty} \int_{-2}^2\int_{-2}^2 H(x,y)\Psi_n(x,y)dxdy=\frac{1}{4\pi}\int_{-2}^2\int_{-2}^2
\frac{(4-xy)H(x,y)}{\sqrt{4-x^2}\sqrt{4-y^2}}dxdy
\end{eqnarray*}
\end{lemma}

Notice that this lemma already implies Theorem \ref{clt:gue} in the case of Lipshitz functions.

\section{Proof of Theorem \ref{clt:gue}}

\subsection{Proof of Necessity}

We will first prove the necessity implication of Theorem \ref{clt:gue}. Sufficiency is dealt with in later sections.

\begin{lemma}
Let $f$ be a continous function on [-2,2] and suppose:
\begin{eqnarray*}
\int \int \left(\frac{f(x)-f(y)}{x-y}\right)^2\frac{4-xy}{\sqrt{4-x^2}\sqrt{4-y^2}}dxdy=\infty
\end{eqnarray*}
Then $\emph{Var}(N_{n}[f])$ diverges in the limit.
\end{lemma}

\textbf{Proof:}

Let $f_j$ be a Lipshitz approximation to $f$, converging uniformly, and let $S$ denote the square $|x|,|y|\leq 2$. Fix $\sigma >0$. By the non-negativity of the integrand:

\begin{eqnarray*}
n\int\int_{x,y\in S} (f(x)-f(y))^2K_n(\sqrt{n}x,\sqrt{n}y)^2dxdy
\\
\geq
n\int\int_{x,y\in S;|x-y|\geq\sigma}(f_j(x)-f_j(y))^2K_n(\sqrt{n}x,\sqrt{n}y)^2dxdy
\\+
n\int\int_{x,y\in S;|x-y|\geq\sigma} ((f(x)-f(y))^2-(f_j(x)-f_j(y))^2)K_n(\sqrt{n}x,\sqrt{n}y)^2dxdy
\end{eqnarray*}

Now, by the uniformness of the approximation we have:
\begin{eqnarray*}
n\int\int_{x,y\in S;|x-y|\geq\sigma} ((f(x)-f(y))^2-(f_j(x)-f_j(y))^2)K_n(\sqrt{n}x,\sqrt{n}y)^2dxdy\\
\geq -\epsilon_j\int\int_{x,y\in S;|x-y|\geq\sigma}nK_n(\sqrt{n}x,\sqrt{n}y)^2dxdy
\end{eqnarray*}

Since $1_{|x-y|\geq \sigma}\left(\frac{1}{x-y}\right)^2$ is a bounded function, we may use Lemma \ref{conv} to obtain:
\begin{eqnarray*}
\lim_{n \to \infty} n\int\int_{x,y\in S;|x-y|\geq\sigma}K_n(\sqrt{n}x,\sqrt{n}y)^2dxdy\\
=\frac{1}{4\pi^2} \int\int_{x,y\in S;|x-y|\geq\sigma}\frac{1}{(x-y)^2} \frac{4-xy}{\sqrt{4-x^2}\sqrt{4-y^2}}dxdy
\end{eqnarray*}

The second integral is finite, so we have the following estimate unformly in $n$:
\begin{eqnarray*}
n\int\int_{x,y\in S;|x-y|\geq\sigma}K_n(\sqrt{n}x,\sqrt{n}y)^2dxdy
\leq C(\sigma)
\end{eqnarray*}

Consequently, we obtain for any fixed $j$ and $\sigma$:
\begin{eqnarray*}
\lim_{n\to \infty}n\int_{-2}^{2}\int_{-2}^{2} (f(x)-f(y))^2K_n(\sqrt{n}x,\sqrt{n}y)^2dxdy
\geq -C(\sigma)\epsilon_j\\
+\limsup_{n\to\infty}n\int\int_{x,y\in S;|x-y|\geq\sigma}(f_j(x)-f_j(y))^2K_n(\sqrt{n}x,\sqrt{n}y)^2dxdy
\end{eqnarray*}

Applying Lemma \ref{conv} again, for any fixed $j$ and $\sigma$:
\begin{eqnarray*}
\lim_{n\to \infty}n\int_{-2}^{2}\int_{-2}^{2} (f(x)-f(y))^2K_n(\sqrt{n}x,\sqrt{n}y)^2dxdy
\geq -C\epsilon_j\\
+
\frac{1}{4\pi^2}\int\int_{x,y\in S; |x-y|\geq \sigma} \left(\frac{f_j(x)-f_j(y)}{x-y}\right)^2\frac{4-xy}{\sqrt{4-x^2}\sqrt{4-y^2}}dxdy
\end{eqnarray*}

With $\sigma$ still fixed, we take the supremum with respect to $j$:
\begin{eqnarray*}
\lim_{n\to \infty}n\int_{-2}^{2}\int_{-2}^{2} (f(x)-f(y))^2K_n(\sqrt{n}x,\sqrt{n}y)^2dxdy
\\
\geq
\frac{1}{4\pi^2}\int\int_{x,y\in S; |x-y|\geq \sigma} \left(\frac{f(x)-f(y)}{x-y}\right)^2\frac{4-xy}{\sqrt{4-x^2}\sqrt{4-y^2}}dxdy
\end{eqnarray*}

Finally, we may use the monotone convergence theorem to send $\sigma$ to zero to obtain:

\begin{eqnarray*}
\lim_{n\to \infty}n\int_{-2}^{2}\int_{-2}^{2} (f(x)-f(y))^2K_n(\sqrt{n}x,\sqrt{n}y)^2dxdy
\\
\geq
\frac{1}{4\pi^2}\int_{-2}^{2}\int_{-2}^{2} \left(\frac{f(x)-f(y)}{x-y}\right)^2\frac{4-xy}{\sqrt{4-x^2}\sqrt{4-y^2}}dxdy
\end{eqnarray*}

This is the desired result. $\P$

\subsection{Sufficiency in the Bulk}

Now that we have necessity, it remains to prove sufficiency. First, we will show sufficiency for the convergence of the variance.
Once this is established, it will not be hard to extend this to the convergence of the linear statistic $N_n[f]$ itself.

We will first consider the case where $f$ is supported away from the spectral edge:

\begin{lemma}
If $f$ is compactly supported on $[-2+\delta,2-\delta]$ for some $\delta>0$, then we have the following equality if the right hand side is finite:
\begin{equation}
\lim_{n\to\infty} Var(N_{n}(f))=\frac{1}{4\pi^2}\int \int \left(\frac{f(x)-f(y)}{x-y}\right)^2\frac{4-xy}{\sqrt{4-x^2}\sqrt{4-y^2}}dxdy
\end{equation}
If the right hand side diverges, so does the left.
\end{lemma}

\textbf{Proof:}

Since we already shown this to be true if the right hand side diverges, we may assume that it is finite and consequently that $f$ is continuous. For convenience, define:
\begin{eqnarray*}
F(x,y)=\left(\frac{f(x)-f(y)}{x-y}\right)^2
\end{eqnarray*}

 Recall that we have the following formula for the variance:

\begin{eqnarray*}
\mbox{Var}(N_n[f])=n\int_{-2,2}\int_{-2,2} F(x,y)
\Big(
\psi_n(\sqrt{n}x)^2\psi_n(\sqrt{n}y)^2
\\
-
\psi_n(\sqrt{n}x)\psi_{n-1}(\sqrt{n}x)\psi_n(\sqrt{n}y)\psi_{n-1}(\sqrt{n}y)
\Big)
 dxdy
\end{eqnarray*}

For now, let's focus on the first term. Fix some $y\in[-2,2]$, and for some sufficiently small $\gamma >0$ write:
\begin{eqnarray*}
\sqrt{n}\int_{-2}^{2}F(x,y)\psi_{n}(\sqrt{n}x)^2dx=
\sqrt{n}\int_{B(y,\gamma)}F(x,y)\psi_{n}(\sqrt{n}x)^2dx
\\
+\sqrt{n}\int_{B(y,\gamma)^c}F(x,y)\psi_{n}(\sqrt{n}x)^2dx
\end{eqnarray*}

Here, $B(y,\gamma)$ denotes the ball of radius $\gamma$ about $y$. Notice that on $B(y,\gamma)^c$ the function $F(x,y)$ is bounded, and consequently by Lemma \ref{conv}:
\begin{eqnarray*}
\lim_{n \to \infty}
\sqrt{n}\int_{B(y,\gamma)^c}F(x,y)\psi_{n}(\sqrt{n}x)^2dx=
\int_{B(y,\gamma)^c}F(x,y)\frac{1}{\pi\sqrt{4-x^2}}dx
\end{eqnarray*}

By the proof of Lemma \ref{conv}, we can take this convergence to be uniform, at least for $y$ away from the edge. On the other hand, since $\sqrt{n}\psi_n(\sqrt{n}x)^2$ is bounded in the bulk, we have:
\begin{eqnarray*}
\sqrt{n}\int_{B(y,\gamma)}F(x,y)\psi_{n}(\sqrt{n}x)^2dx\leq
C(y)\int_{B(y,\gamma)}F(x,y) dx
\end{eqnarray*}

Since $F$ is integrable we have that, for any fixed $\epsilon>0$, we can choose $n$ large enough so that:
\begin{eqnarray*}
\left| \sqrt{n}\int_{-2}^{2}F(x,y)\psi_{n}(\sqrt{n}x)^2dx-\int \frac{F(x,y)}{\pi\sqrt{4-x^2}}dx\right|\leq \epsilon/4
\end{eqnarray*}

Notice that the constant $C(y)$ is uniformly bounded for $y$ in the interval $\Omega_{\alpha}=[-2+\alpha, 2-\alpha]$, so in this interval we can in fact choose $n(\alpha)$ large enough that:
\begin{eqnarray*}
\Big|n\int_{\Omega(\alpha)}\left(\frac{1}{\pi}\int_{-2}^{2}\left(\frac{f(x)-f(y)}{x-y}\right)^2
\frac{1}{\sqrt{4-x^2}}dx\right)\psi_n(\sqrt{n}y)^2 dy
\\
-\sqrt{n}\int_{\Omega(\alpha)}\left(\int_{-2}^{2}\frac{F(x,y)}{\pi\sqrt{4-x^2}}dx\right)
\psi_n(\sqrt{n}y)^2dy\Big|\leq \epsilon/4
\end{eqnarray*}

We can choose $\alpha$ so that $\Omega(\alpha)$ compactly contains the support of $f$, so that:
\begin{eqnarray*}
n\int_{\Omega(\alpha)^c}\left(\frac{1}{\pi}\int_{-2}^{2}\left(\frac{f(x)-f(y)}{x-y}\right)^2
\frac{1}{\sqrt{4-x^2}}dx\right)\psi_n(\sqrt{n}y)^2 dy\\=
n\int_{\Omega(\alpha)^c}\left(\int_{-2}^{2}\left(\frac{f(x)}{x-y}\right)^2\psi_n(\sqrt{n}x)^2 dx\right) \psi_n(\sqrt{n}y)^2dy
\end{eqnarray*}

For $y\in\Omega(\alpha)^c$ we may bound:
\begin{eqnarray*}
\left(\frac{f(x)}{x-y}\right)^2\leq C(\alpha)
\end{eqnarray*}

This is because, for $y$ in this region, if $f(x)$ does not vanish then $(x-y)^{-1}$ cannot be too large (with the exact constant depending on how much $\Omega(\alpha)$ overshoots the support of $f$). By Lemma \ref{conv}, then:

\begin{eqnarray*}
\limsup_{n\to\infty}n\int_{\Omega(\alpha)^c}\left(\frac{1}{\pi}\int_{-2}^{2}\left(\frac{f(x)-f(y)}{x-y}\right)^2
\frac{1}{\sqrt{4-x^2}}dx\right)\psi_n(\sqrt{n}y)^2 dy
\\
\leq
\frac{1}{4\pi^2}\int_{-2}^{2}\int_{-2}^{2} 1_{\Omega(\alpha)}(y)\left(\frac{f(x)}{x-y}\right)^2
\frac{4-xy}{\sqrt{4-x^2}\sqrt{4-y^2}}dxdy
\end{eqnarray*}

This quantity vanishes with $\alpha$ (by the assumption of the finiteness of the limiting variance). Therefore, we may choose $\alpha$ and $n(\alpha)$ so that:
\begin{eqnarray*}
\left|n\int_{\Omega(\alpha)^c}\left(\frac{1}{\pi}\int_{-2}^{2}\left(\frac{f(x)-f(y)}{x-y}\right)^2
\frac{1}{\sqrt{4-x^2}}dx\right)\psi_n(\sqrt{n}y)^2 dy\right|\leq \epsilon/8
\end{eqnarray*}

By a similar argument, we can simulataneously arrange $\alpha$ and $n$ so that:
\begin{eqnarray*}
|\sqrt{n}\int_{\Omega(\alpha)^c}\left(\int_{-2}^{2}\frac{F(x,y)}{\pi\sqrt{4-x^2}}dx\right)
\psi_n(\sqrt{n}y)^2dy|\leq \epsilon/8
\end{eqnarray*}

Consequently, for large enough $n$:
\begin{eqnarray*}
\Big|n\int_{-2}^{2}\left(\frac{1}{\pi}\int_{-2}^{2}\left(\frac{f(x)-f(y)}{x-y}\right)^2
\frac{1}{\sqrt{4-x^2}}dx\right)\psi_n(\sqrt{n}y)^2 dy
\\
-\sqrt{n}\int_{-2}^{2}\left(\int_{-2}^{2}\frac{F(x,y)}{\pi\sqrt{4-x^2}}dx\right)
\psi_n(\sqrt{n}y)^2dy\Big|\leq \epsilon/2
\end{eqnarray*}

We now need to essentially repeat these steps to handle the integration in $y$. For convenience, define:

\begin{eqnarray*}
H(y)=\int_{-2}^{2}\frac{F(x,y)}{\pi\sqrt{4-x^2}}dx
=\frac{1}{\pi}\int_{-2}^{2}\left(\frac{f(x)-f(y)}{x-y}\right)^2
\frac{1}{\sqrt{4-x^2}}dx
\end{eqnarray*}

Let $\Omega(\sigma)$ denote the interval $[-2+\sigma, 2-\sigma]$ for $\sigma$ small . As with $\Omega(\alpha)$, we take $\sigma$ small enough so that $\Omega(\sigma)$ compactly contains the support of $f$), and consequently, for $y\in \Omega(\sigma)^c$ we have the estimate $H(y)\leq C(\sigma)$, which combined with lemma 4 provides:
\begin{eqnarray*}
\lim_{n\to\infty}\sqrt{n}\int_{\Omega(\sigma)^c}H(y)\psi_n(\sqrt{n}y)^2dy
=\frac{1}{\pi}\int_{\Omega(\sigma)^c} H(y)
\frac{1}{\sqrt{4-y^2}}dy
\end{eqnarray*}

We concentrate then on the region $\Omega(\sigma)$. Since $F(x,y)$ is integrable, $H(y)$ is integrable. By the Chebycheff inequality, we may divide $\Omega(\sigma)$ into a region on which $H(y)$ is bounded, and an arbitrarily small region on which $H(y)$ is not. The contribution of the arbitrarily small region is again arbitrarily small (since $\psi_n(\sqrt{n}y)^2$ is uniformly bounded in the bulk and $H(y)$ is integrable, and therefore must have small $L^1$ norm on small sets).

Lemma \ref{conv} takes care of the region on which $H(y)$ is bounded, so we may set the error to be no larger than $\epsilon/2$,  so that we can arrange:

\begin{eqnarray*}
\Big| \sqrt{n}\int_{-2}^{2}H(y)\psi_n(\sqrt{n}y)^2dy-\frac{1}{\pi}\int_{-2}^{2}
H(y)\frac{1}{\sqrt{4-y^2}}dy\Big|\leq \epsilon/2
\end{eqnarray*}

Combining with out previous estimate, we get:
\begin{eqnarray*}
\Big|n\int_{-2}^{2}\left(\frac{1}{\pi}\int_{-2}^{2}\left(\frac{f(x)-f(y)}{x-y}\right)^2
\frac{1}{\sqrt{4-x^2}}dx\right)\psi_n(\sqrt{n}y)^2 dy
\\
-\frac{1}{\pi^2}\int_{-2}^{2}\int_{-2}^{2}\left(\frac{f(x)-f(y)}{x-y}\right)^2\frac{1}{\sqrt{4-y^2}\sqrt{4-x^2}}dxdy\Big|\leq \epsilon
\end{eqnarray*}

This handles the first term in the expansion for the variance. The second is handled analagously, proving the result. $\P$

\subsection{Estimates at the Edge, And Convergence of Higher Cumulants}

The next order of business is to extend our work to the edge of spectral support.

Our main tool will be an aymptotic estimate for Hermite polynomials that follows from a paper of Erdelyi \cite{Erd}, who considers the following type of differential equations:

\begin{eqnarray*}
y^{\prime \prime}+[\lambda^2p(x)+r(x,\lambda)]y=0
\end{eqnarray*}

In \cite{Erd} an asymptotic expansion for solutions of these equations is obtained, under some conditions, and applied to several special cases. One of these is the case of Hermite polynomials, and the results obtained imply:

\begin{lemma}
Define $\phi(x)$ by the following formula:
\begin{eqnarray*}
\phi(x)=\frac{3}{2}\left(
\frac{1}{2}x(1-x^2)^{1/2}-\frac{1}{2}\cos^{-1}(x)
\right)^{2/3}
\end{eqnarray*}

The following estimate for Hermite oscillator functions holds:

\begin{eqnarray*}
\label{upbound}
|\psi_n(\sqrt{n}x)|\leq \frac{C}{n^{1/4}}
\left|\phi^{\prime}\left(\frac{\sqrt{n}x}{\sqrt{2n+1}\sqrt{2}}\right)\right|^{-1/2}
\end{eqnarray*}

The constant $C$ is independent of $n$.
\end{lemma}

\textbf{Proof:}

Formula (6.11) in \cite{Erd} gives the following heuristic approximation, valid for $|x|<1$:

\begin{eqnarray*}
2^{n/2}h_n(\sqrt{2}\sqrt{2n+1}x)\approx \sqrt{2\pi}(2n+1)^{n/2+1/6}e^{\frac{2n+1}{2}\times(x^2-\frac{1}{2})}
\times\\
|\phi^{\prime}|^{-1/2}\mbox{Ai}(-(2n+1)^{2/3}\phi)\left(1+O\left(\frac{1}{n(1+x^2)}\right)\right)
\end{eqnarray*}

This formula is not precisely rigorous because of problems by the zeros of the Airy function, but it turns out will suffice for our purposes. This is because we only seek to employ it in the capacity of an upper bound, so in principle we should be able to get away with replacing the Airy function $\mbox{Ai}(z)$ with $z^{1/4}$. In fact, this is exactly what happens.

To be precise, the discussion preceeding (6.11) in \cite{Erd} provides a fully rigorous version of the above approximation:

\begin{eqnarray*}
2^{n/2}h_n(\sqrt{2}\sqrt{2n+1}x)= \sqrt{2\pi}(2n+1)^{n/2+1/6}e^{\frac{2n+1}{2}\times(x^2-\frac{1}{2})}y_0(x)
\end{eqnarray*}

By (2.15) in  \cite{Erd}, for $x\in[0,\infty]$ the function $y_0(x)$ uniformly satisfies:

\begin{eqnarray*}
y_0(x)=Y_0(x)+O\left(\tilde{Y}_0(x)\frac{1}{n(1+x^2)}\right)
\end{eqnarray*}

The functions $Y_0(x),\tilde{Y}_0(x)$ themselves uniformly satisfy (by equations (2.10)-(2.11), (3.1)-(3.4) in \cite{Erd}) the inequality:
\begin{eqnarray*}
\max\left[|Y_0(x)|,|\tilde{Y}_0(x)|\right]\leq C|\phi^{\prime}|^{-1/2}|N^{2/3}\phi|^{-1/4}
\end{eqnarray*}

This allows us to turn the heuristic approximation into a rigorous one. For $|x|<1$:

\begin{eqnarray*}
|2^{n/2}h_n(\sqrt{2}\sqrt{2n+1}x)|\leq C (2n+1)^{n/2+1/6}e^{\frac{2n+1}{2}\times(x^2-\frac{1}{2})}
|\phi^{\prime}|^{-1/2}|(2n+1)^{2/3}\phi|^{-1/4}
\end{eqnarray*}

Rearranging:
\begin{eqnarray*}
\left|\left(\frac{1/n}{1+1/2n}\right)^{n/2}h_n(\sqrt{2}\sqrt{2n+1}x)\right|\leq C e^{(n+1/2)(x^2-\frac{1}{2})}
|\phi^{\prime}|^{-1/2}|\phi|^{-1/4}
\end{eqnarray*}

Taking the limit of $(1+1/2n)^{n/2}$, replacing $x$ with $\sqrt{n}x/(\sqrt{2}\sqrt{2n+1})$ for $|x|< 2$ (which still has magnitude less than one):

\begin{eqnarray*}
\left|\left(\frac{1}{n}\right)^{n/2}h_n(\sqrt{n}x)\right|\leq C e^{-n/2}e^{n^2x^2/(4n+2)}
|\phi^{\prime}|^{-1/2}|\phi|^{-1/4}
\end{eqnarray*}

Substituting in the defintion  $\psi_n(x)=e^{-x^2/4}h_n(x)(\sqrt{n!}(2\pi)^{1/4})^{-1}$:

\begin{eqnarray*}
\left|\left(\frac{1}{n}\right)^{n/2}e^{nx^2/4}\psi_n(\sqrt{n}x)\sqrt{n!}\right|\leq C e^{-n/2}e^{n^2x^2/(4n+2)}
|\phi^{'}|^{-1/2}|\phi|^{-1/4}
\end{eqnarray*}

By the Stirling approximation:

\begin{eqnarray*}
\left|e^{nx^2/4}h_n(\sqrt{n}x)n^{1/4}\right|\leq C e^{n^2x^2/(4n+2)}
|\phi^{\prime}|^{-1/2}|\phi|^{-1/4}
\end{eqnarray*}

Applying the limit $n/(4+2/n)-n/4\to1/8$ as $n$, and also $|\phi|=O(1)$, we have:

\begin{eqnarray*}
|n^{1/4}\psi_n(\sqrt{n}x)|\leq C
\left|\phi^{\prime}\left(\frac{\sqrt{n}x}{\sqrt{2n+1}\sqrt{2}}\right)\right|^{-1/2}
\end{eqnarray*}

This is as we sought. $\P$.

We are now ready to study the spectral edge. Define the following region:
\begin{eqnarray*}
\Omega_\delta=\{x\in[2-\delta,2]\cup[-2,-2+\delta]\}\cup
\{y\in[2-\delta,2]\cup[-2,-2+\delta]\}
\end{eqnarray*}

We will show that this region does not contribute much to the limiting variance . Specifically, we may control contribution made by a neighborhood of the edge by controling the size of the neighborhood.

\begin{theorem} Let $\delta>0$. There exists a constant $C$, which does not depend on $\delta$, such that for any continuous function $f$ such that $V_{GUE}[f]< \infty$, the following estimate holds:
\begin{eqnarray*}
\Bigg|n\int\int_{\Omega_\delta}\left(\frac{f(x)-f(y)}{x-y}\right)^2
\Big(\psi_n(\sqrt{n}x)^2\psi_n(\sqrt{n}y)^2\\
+\psi_n(\sqrt{n}x)\psi_{n-1}(\sqrt{n}x)\psi_n(\sqrt{n}y)\psi_{n-1}(\sqrt{n}y)\Big)
dxdy\Bigg|\\
\leq C\int\int_{\Omega_\delta}\left(\frac{f(x)-f(y)}{x-y}\right)^2\frac{4-xy}{\sqrt{4-x^2}\sqrt{4-y^2}}dxdy
\end{eqnarray*}
Note that if $\delta$ is chosen to be small enough, the right hand side can be made less than any $\epsilon>0$.
\end{theorem}

\textbf{Proof:}

In what follows, we will focus on the first summand, featuring the term $\psi_n(\sqrt{n}x)^2\psi_n(\sqrt{n}y)^2$. The argument given will not change if (for instance) we change a $\psi_n(\sqrt{n}x)^2$ to $\psi_{n-1}(\sqrt{n}x)^2$, so by the Cauchy-Young inequality controlling the first term is sufficient to prove the full result. We also assume that $x\approx 2$ and $y \approx 2$ (since Hermite polynomials are either even or odd, it suffices to consider this case).

Start by taking a simple derivative:
\begin{eqnarray*}
\phi^{\prime}(x)=\left(\frac{x}{2}(1-x^2)^{1/2}-\frac{1}{2}\cos^{-1}(x)\right)^{-1/3}(1-x^2)^{1/2}
\end{eqnarray*}

Consequently, the upper bound in Lemma \ref{upbound} can be written:
\begin{eqnarray*}
\left|\phi^{\prime}\left(\frac{\sqrt{n}x}{\sqrt{2n+1}\sqrt{2}}\right)\right|^{-1/2}=
\frac{\left(\frac{x}{\sqrt{4+2/n}}(1-\frac{x^2}{4+2/n})^{1/2}-\cos^{-1}(\frac{x}{\sqrt{4+2/n}})
\right)^{1/6}}
{2^6\times(1-\frac{x^2}{4+2/n})^{1/4}}
\end{eqnarray*}

Suppose now that we had the ratio bound $|R(x,y,n)^{1/3}|\leq C$, where $R(x,y,n)$ denotes the following expression:

\begin{eqnarray*}
\frac{\left(\frac{x}{4+2/n}(4-x^2+2/n)^{1/2}-\cos^{-1}(\frac{x}{\sqrt{4+2/n}})
\right)\left(\frac{y}{4+2/n}(4-y^2+2/n)^{1/2}-\cos^{-1}(\frac{y}{\sqrt{4+2/n}})
\right)}
{(4-xy+\frac{2}{n})^3}
\end{eqnarray*}

This and the previous asymptotic expansion would allow us to write:

\begin{eqnarray*}
|n^{1/2}\psi_n(\sqrt{n}x)^2\times n^{1/2}\psi_n(\sqrt{n}y)^2|
\leq
C_1|R(x,y,n)^{1/3}|
\frac{4+\frac{2}{n}-xy}{\sqrt{4-y^2+2/n}\sqrt{4-x^2+2/n}}\\
\leq C_2
\frac{4-xy}{\sqrt{4-y^2}\sqrt{4-x^2}}
\end{eqnarray*}

This in turn gives us the sought-after result, so it only remains to prove the claimed ratio bound. For small $x$, the inverse of cosine may be expanded as follows:
\begin{eqnarray*}
\cos^{-1}(1-x)=\sqrt{2x}+\frac{1}{24}(2x)^{3/2}+O(x^{2})
\end{eqnarray*}

In our case, this becomes:
\begin{eqnarray*}
\cos^{-1}\left(\frac{x}{\sqrt{4+\frac{2}{n}}}\right)=
\sqrt{2-\frac{2x}{\sqrt{4+\frac{2}{n}}}}+\frac{1}{24}\left(2-\frac{2x}{\sqrt{4+\frac{2}{n}}}\right)^{3/2}
+O\left[\left(1-\frac{x}{\sqrt{4+\frac{2}{n}}}\right)^{2}\right]
\end{eqnarray*}

We can then write $R(x,y,n)$ as:
\begin{eqnarray*}
\left(\frac{x}{4+\frac{2}{n}}(4+\frac{2}{n}-x^2)^{1/2}
-\frac{\sqrt{2}}{(4+\frac{2}{n})^{1/4}}\sqrt{\sqrt{4+\frac{2}{n}}-x}-
\frac{1}{24}\left(2-\frac{2x}{\sqrt{4+\frac{2}{n}}}\right)^{3/2}
+O\left[\left(1-\frac{x}{\sqrt{4+\frac{2}{n}}}\right)^{2}\right]
\right)\\
\left(\frac{y}{4+\frac{2}{n}}(4+\frac{2}{n}-y^2)^{1/2}
-\frac{\sqrt{2}}{(4+\frac{2}{n})^{1/4}}\sqrt{\sqrt{4+\frac{2}{n}}-y}-
\frac{1}{24}\left(2-\frac{2y}{\sqrt{4+\frac{2}{n}}}\right)^{3/2}
+O\left[\left(1-\frac{y}{\sqrt{4+\frac{2}{n}}}\right)^{2}\right]
\right)\\
\times \left(\frac{1}{4+\frac{2}{n}-xy}\right)^3
\end{eqnarray*}

We will ignore the error terms for now, and concentrate on bounding the other terms in the expression. Writing $(4+\frac{2}{n}-x^2)^{1/2}=\sqrt{(4+\frac{2}{n})^{1/2}-x}\times\sqrt{(4+\frac{2}{n})^{1/2}+x}$, we have (up to error terms):
\begin{eqnarray*}
R(x,y,n)\approx\frac{\sqrt{(4+\frac{2}{n})^{1/2}-x}\sqrt{(4+\frac{2}{n})^{1/2}-y}}{4+\frac{2}{n}-xy}
\times\frac{1}
{(4+\frac{2}{n}-xy)^2}
\\
\times\left(\frac{x}{4+\frac{2}{n}}\sqrt{(4+\frac{2}{n})^{1/2}+x}-
\frac{\sqrt{2}}{(4+2/n)^{1/4}}-\frac{1}{24}\left(2-\frac{2x}{\sqrt{4+\frac{2}{n}}}\right)
\right)\\
\times
\left(\frac{y}{4+\frac{2}{n}}\sqrt{(4+\frac{2}{n})^{1/2}+y}-
\frac{\sqrt{2}}{(4+2/n)^{1/4}}-\frac{1}{24}\left(2-\frac{2y}{\sqrt{4+\frac{2}{n}}}\right)
\right)
\end{eqnarray*}

Using the trivial $x^2+y^2-2xy\geq 0$, we have the following inequality:
\begin{eqnarray*}
0\leq
\frac{\sqrt{(4+\frac{2}{n})-x^2}\sqrt{(4+\frac{2}{n})-y^2}}{4+\frac{2}{n}-xy}\leq 1
\end{eqnarray*}

It remains to control the other factors. We start with the following estimate:

\begin{eqnarray*}
\frac{\left(\frac{x}{4+\frac{2}{n}}\sqrt{(4+\frac{2}{n})^{1/2}+x}-
\frac{\sqrt{2}}{(4+2/n)^{1/4}}\right)
\left(\frac{y}{4+\frac{2}{n}}\sqrt{(4+\frac{2}{n})^{1/2}+y}-
\frac{\sqrt{2}}{(4+2/n)^{1/4}}\right)}
{(4+\frac{2}{n}-xy)^2}\\
\leq
\frac{\left(\frac{x}{4+\frac{2}{n}}\sqrt{(4+\frac{2}{n})^{1/2}+x}-
\frac{\sqrt{2}}{(4+2/n)^{1/4}}\right)}{4+\frac{2}{n}-2x}
\times
\frac{\left(\frac{y}{4+\frac{2}{n}}\sqrt{(4+\frac{2}{n})^{1/2}+y}-
\frac{\sqrt{2}}{(4+2/n)^{1/4}}\right)}{4+\frac{2}{n}-2y}
\end{eqnarray*}

Also, setting $u=\frac{x}{\sqrt{4+(2/n)}}$ (so $1/2 < u <1$):
\begin{eqnarray*}
\left|\frac{\frac{x}{4+2/n}\sqrt{(4+2/n)^{1/2}+x}-\frac{\sqrt{2}}{(4+2/n)^{1/4}}}
{(4+2/n)-2x}\right|\\=
\left|\frac{1}{u(4+(2/n))^{3/4}}
\frac{\sqrt{2}-u\sqrt{1+u}}
{\sqrt{2}-u^{-1}(2+1/n)^{1/2}}\right|\\
\leq \left|\frac{2}{(4+(2/n))^{3/4}} \right|\leq C
\end{eqnarray*}

Together, these imply that in our region of interest (that is, $|x|,|y|\leq 2$):
\begin{eqnarray*}
\left|\frac{\left(\frac{x}{4+\frac{2}{n}}\sqrt{(4+\frac{2}{n})^{1/2}+x}-
\frac{\sqrt{2}}{(4+2/n)^{1/4}}\right)
\left(\frac{y}{4+\frac{2}{n}}\sqrt{(4+\frac{2}{n})^{1/2}+y}-
\frac{\sqrt{2}}{(4+2/n)^{1/4}}\right)}
{(4+\frac{2}{n}-xy)^2}\right|
\leq
C
\end{eqnarray*}

Additionally, we have:

\begin{eqnarray*}
\frac{1}{(4+\frac{2}{n}-xy)^2}\left(2-\frac{2x}{\sqrt{4+2/n}}\right)
\left(2-\frac{2y}{\sqrt{4+2/n}}\right)\\
\leq C \frac{\left(\sqrt{4+2n}-x\right)}{(4+\frac{2}{n}-xy)}
\frac{\left(\sqrt{4+2n}-y\right)}{(4+\frac{2}{n}-xy)}\\
\leq C \left(\frac{\sqrt{4+2n}-x}{\sqrt{4+\frac{2}{n}}-\sqrt{xy}}\right)
\left(\frac{\sqrt{4+2n}-y}{\sqrt{4+\frac{2}{n}}-\sqrt{xy}}\right)\leq C
\end{eqnarray*}

And also:
\begin{eqnarray*}
0\leq \left|\frac{\left(2-\frac{2x}{\sqrt{4+2/n}}\right)
\left(\frac{y}{4+\frac{2}{n}}\sqrt{(4+\frac{2}{n})^{1/2}+y}-
\frac{\sqrt{2}}{(4+2/n)^{1/4}}\right)}
{(4+2/n-xy)^2}\right|\\
 \leq
\left|
\frac{
\left(\frac{y}{4+\frac{2}{n}}\sqrt{(4+\frac{2}{n})^{1/2}+y}-
\frac{\sqrt{2}}{(4+2/n)^{1/4}}\right)}
{4+2/n-2y}\right|
\times\left|\frac{\left(2-\frac{2x}{\sqrt{4+2/n}}\right)}
{4+2/n-2x}\right|
\\
\leq C\left|\frac{1-\frac{x}{\sqrt{4+2/n}}}
{4+2/n-2x}\right|\leq 4C
\end{eqnarray*}

The last equality is due to the positivity of $6+x-4\sqrt{2x}$ on $[1,2]$. Combining all of these estimates, we have now verified the bound:
\begin{eqnarray*}
\Bigg|\frac{1}
{(4+\frac{2}{n}-xy)^2}
\times
\left(\frac{x}{4+\frac{2}{n}}\sqrt{(4+\frac{2}{n})^{1/2}+x}-
\frac{\sqrt{2}}{(4+2/n)^{1/4}}-\frac{1}{24}\left(2-\frac{2x}{\sqrt{4+\frac{2}{n}}}\right)
\right)\\
\times
\left(\frac{y}{4+\frac{2}{n}}\sqrt{(4+\frac{2}{n})^{1/2}+y}-
\frac{\sqrt{2}}{(4+2/n)^{1/4}}-\frac{1}{24}\left(2-\frac{2y}{\sqrt{4+\frac{2}{n}}}\right)
\right)\Bigg|
\leq C
\end{eqnarray*}

The terms which contain factors of the error terms of the $\cos^{-1}$ expansion are controlled by the same estimates, so taking a cubic root finishes the proof $\P$.

We now have the convergence of the variance of the linear statistic to the given limit under optimal conditions. To prove Theorem \ref{clt:gue}, all that is left to do is to show that the convergence of the variance implies the convergence of the linear statistic itself (to a Gaussian limit).  Since the result is already proven for smooth functions, this follows from an approximation argument.

\textbf{Proof of Theorem \ref{clt:gue}:}

Let $\Omega=[-2,2]\times[-2,2]$, and let $f(x)$ be a function on $[-2,2]$ such that $V_{\mbox{GUE}}[f]$ is finite.
Then necessarily $f$ is continuous -- indeed, we have the Sobolev bound:

\begin{eqnarray*}
\int_{-2}^2\int_{-2}^2 \left(\frac{f(x)-f(y)}{x-y}\right)^2dxdy \leq C||f||_{H^{1/2}}
\end{eqnarray*}

Fix $\epsilon >0$. For some small $\delta>0$, let $f_k$ be  a sequence of smooth functions converging to $f$ in the $H^{1/2}[-2+\delta,2-\delta]$ norm.
Since we may take $f_k$ to be mollifications of $f$ without loss of generality, we may take $|f_k(x)|<||f||_{\infty[-2,2]}$.  Then we have (for $k$ large enough):
\begin{eqnarray*}
\int_{-2}^2\int_{-2}^2 \left(\frac{(f(x)-f_k(x))-(f(y)-f_k(y))}{x-y}\right)^2\frac{4-xy}{\sqrt{4-x^2}
\sqrt{4-y^2}}dxdy < \epsilon^3
\end{eqnarray*}

Indeed, fix small (smaller than $\delta$) $\delta_1 > \delta_2 > 0$. When $x$ and $y$ both stay away from the edge (say, $\min[|x-2|,|y-2|]>\delta_2$) the smallness of this quantity follows directly from $H^{1/2}$ approximation. When $x$ and $y$
are both close to the edge (for instance, say $\max[|x-2|,|y-2|<\delta_1]$ holds) $f_k(x)$ and $f_k(y)$ both vanish and smallness follows from the fact that $V_{GUE}$ is finite (and therefore the integral must be small on small sets).

When either $x$ is close to the edge and the $y$ isn't (so $|x-2|<\delta_2,|y-2|>\delta_1$) or vice versa, we can bound $|x-y|^{-1}$ (since there is some space between $\delta_1$ and $\delta_2$), and the result then follows
from the fact that the function $\frac{4-xy}{\sqrt{4-x^2}\sqrt{4-y^2}}$ is integrable. Therefore, we may write:
\begin{eqnarray*}
\lim_{n\to\infty}\mbox{Var}[N_n(f-f_k)]\leq \epsilon^3
\end{eqnarray*}

Using the linearity of the linear statistic, the result follows easily from approximation (using either characteristic functions or by metrizing convergence in distribution, for instance by employing the Ky Fan metric). $\P$

\newpage

\section{Appendix}

\subsection{Proof of Lemma \ref{conv}}

In this section we prove Lemma \ref{conv}. We will first prove the first two equalities asserted in the lemma, and then the third (which is a conseqeuece of the first two).

\begin{lemma}
The first and second equations of Lemma \ref{conv} hold for all functions in $L^{\infty}[-2,2]$.
\end{lemma}

\textbf{Proof:}
First assume that $\phi$ is smooth function with compact support in $(-2,2)$. By the known asymptotics of $K^n$ (see for instance \cite{Tao}):
\begin{equation}
\lim\frac{1}{\sqrt{n}}\int \phi(x)K^{n}(\sqrt{n}x,\sqrt{n}x)dx=\frac{1}{2\pi}\int\phi(x)\sqrt{4-x^2}dx
\end{equation}

By the Christoffel-Darboux formula, the identity $\psi^{\prime}_n=\frac{-x}{2}\psi_n+\sqrt{n}\psi_{n-1}$, and integration by parts, we obtain:

\begin{eqnarray*}
\int \phi(x)K^{n}(\sqrt{n}x,\sqrt{n}x)dx
=\frac{1}{2n}\int\phi^{''}(x)\psi_n(\sqrt{n}x)^2dx
\\
-\int\phi(x)[2\psi_n(\sqrt{n}x)\psi^{''}_n(\sqrt{n}x)+\frac{1}{2}\psi_n(\sqrt{n}x)^2]dx
\end{eqnarray*}

Now we apply the formula $\psi^{''}_n=-(n+\frac{1}{2}-\frac{x^2}{4})\psi_n$:
\begin{eqnarray*}
\frac{1}{\sqrt{n}}\int \phi(x)K^{n}(\sqrt{n}x,\sqrt{n}x)dx=
\frac{1}{2\sqrt{n}n}\int\phi^{''}(x)\psi_n(\sqrt{n}x)^2dx
\\
+\sqrt{n}\int \phi(x)[2-\frac{x^2}{2}]\psi_n(\sqrt{n}x)^2dx+
\frac{1}{\sqrt{n}}\int \psi_n(\sqrt{n}x)^2dx
\end{eqnarray*}

Passing to the limit:

\begin{equation*}
\lim_{n\to\infty} \frac{1}{2}\sqrt{n}\int \phi(x)(4-x^2)\psi_n(\sqrt{n}x)^2dx=\frac{1}{2\pi}\int \phi(x)\sqrt{4-x^2}dx
\end{equation*}

Writing $\frac{f(x)}{4-x^2}=\phi(x)$ (which is permissible since we have compact support in $(-2,2)$) we obtain:

\begin{equation*}
\lim_{n\to\infty} \sqrt{n}\int f(x)\psi_n(\sqrt{n}x)^2dx=\frac{1}{\pi}\int \frac{f(x)}{\sqrt{4-x^2}}dx
\end{equation*}

If $f$ is a nonnegative function supported in $[a,b]\in(-2,2)$ and continuous on $[a,b]$, then the value of $\limsup\int f(x)\sqrt{n}\psi_n(\sqrt{n}x)^2dx$ is bounded above by $\frac{1}{\pi}\int \frac{F(x)}{\sqrt{4-x^2}}dx$, where $F$ is any extension of $f$ to a smooth, compactly supported nonnegative function on $(-2,2)$ such that than $F(x)\geq f(x)$ holds everywhere. Taking the infinimum of the values, we see that:
\begin{equation*}
\frac{1}{\pi}\int \frac{f(x)}{\sqrt{4-x^2}}dx\geq \limsup \int f(x)\sqrt{n}\psi_n(\sqrt{n}x)^2dx
\end{equation*}

On the other hand, by multiplying $f$ by an arbitrary cut-off function and arguing similarly we obtain the opposite bound:

\begin{equation*}
\frac{1}{\pi}\int \frac{f(x)}{\sqrt{4-x^2}}dx\leq \liminf \int f(x)\sqrt{n}\psi_n(\sqrt{n}x)^2dx
\end{equation*}

Thus, we have enlarged our space of test functions to include, for instance, simple functions (sums of multiples of indicator functions) supported on compactly on $(-2,2)$.

For any $f\in L^\infty[a,b]$, where $[a,b]$ is still a proper compact subset of $(-2,2)$, let $f_j$ be a sequence of simple functions which converges to $f$ in the $L^{\infty}$ norm.  By the triangle inequality:

\begin{eqnarray*}
\left|\int_{a}^{b} f(x)\sqrt{n}\psi_n(\sqrt{n}x)^2dx-\frac{1}{\pi}\int_{a}^{b} f(x)\frac{1}{\sqrt{4-x^2}}dx\right|
\leq\int_{a}^{b} |f(x)-f_j(x)|\sqrt{n}\psi_n(\sqrt{n}x)^2dx\\+
\left|\int_{a}^{b} f_j(x)\sqrt{n}\psi_n(\sqrt{n}x)^2dx-\frac{1}{\pi}\int_{a}^{b} f_j(x)\frac{1}{\sqrt{4-x^2}}dx\right|
+\frac{1}{\pi}\int_{a}^{b}|f_j(x)-f(x)|\frac{1}{\sqrt{4-x^2}}dx
\end{eqnarray*}

Choosing $j$ and $n(j)$ large enough forces the right hand side to be less than $\epsilon$, so we now have the result for essentially bounded function with support in the bulk. Finally, we lose the assumption of compact support. Take $f\in L_{\infty}$ and let $f_\epsilon$ be the restriction of $f$ to the interval $(-2+t,2-t)$, and $f^c_t=f-f_t$. We may expand:
\begin{equation*}
\int f(x)\sqrt{n}\psi^2_n(\sqrt{n}x)dx=\sqrt{n}\int f_\epsilon(x) \psi^2_n(\sqrt{n}x)dx +
\sqrt{n}\int f^c_t(x)\psi^2_n(\sqrt{n}x)dx
\end{equation*}

Since we already have the result for indicator functions and $|f^c_t(x)|\leq C$, we can make the second term on the right hand side arbitrarily small (in the limit as $n$ goes to infinity) by taking $t$ small. Of course, the result holds for $f_t$ because this function has compact support in $(-2,2)$, so we have:
\begin{equation*}
\limsup\int f(x)\sqrt{n}\psi_n(\sqrt{n}x)^2dx\leq\frac{1}{\pi}\int f(x)\frac{1}{\sqrt{4-x}}dx
\end{equation*}

On the other hand, by considering truncations of $f$:
\begin{equation*}
\liminf\int f(x)\sqrt{n}\psi_n(\sqrt{n}x)^2dx\geq\frac{1}{\pi}\int f(x)\frac{1}{\sqrt{4-x}}dx
\end{equation*}

The limits therefore coincide, and we conclude the proof of the first equation in Lemma \ref{conv}.

The proof of the second part of the lemma is very similar. Suppose $f$ is smooth and compactly supported. By the identity $\psi^{'}_n=\frac{-x}{2}\psi_n+\sqrt{n}\psi_{n-1}$ and integration by parts:

\begin{equation*}
\sqrt{n}\int f(x)\psi_n(\sqrt{n}x)\psi_{n-1}(\sqrt{n}x)dx=
\int (\sqrt{n}\frac{x}{2}f(x)-\frac{f'(x)}{2\sqrt{n}})\psi_n(\sqrt{n}x)^2dx
\end{equation*}

Therefore:

\begin{equation*}
\lim \sqrt{n}\int f(x)\psi_n(\sqrt{n}x)\psi_{n-1}(\sqrt{n})dx=\frac{1}{\pi}\int \frac{xf(x)}{2\sqrt{4-x^2}}dx
\end{equation*}

By smooth approximation and the Holder inequality, we may extend this to simple functions, which, just as before, gives us the result for compactly supported $L^\infty$ functions. Finally, we drop the assumption of compact support: take $f\in L^{\infty}$ and let $f_t$ be the restriction of $f$ to the interval $(-2+t,2-t)$, and $f^c_t=f-f_t$. Then:
\begin{eqnarray*}
\int |f(x)-f_{t}(x)|\sqrt{n}\psi_n(\sqrt{n}x)\psi_{n-1}(\sqrt{n}x)dx\\
\leq  ||f||_{\infty}\left(\int_{-2+t}^{2-t}\sqrt{n}\psi^2_n(\sqrt{n}x)dx\right)^2
\left(\int_{-2+t}^{2-t}\sqrt{n}\psi^2_{n-1}(\sqrt{n}x)dx\right)^2
\end{eqnarray*}

Just as before, an $\frac{\epsilon}{3}$ argument finishes the proof. $\P$

Finally, we combine these two formulas to prove:

\begin{lemma}
If $H(x,y)\in L^{\infty}[\mathcal{R}\times \mathcal{R}]$, then:
\begin{eqnarray*}
\lim_{n\to\infty} \int_{-2}^2\int_{-2}^2 H(x,y)\Psi_n(x,y)dxdy=\frac{1}{4\pi}\int_{-2}^2\int_{-2}^2
H(x,y)\frac{4-xy}{\sqrt{4-x^2}\sqrt{4-y^2}}dxdy
\end{eqnarray*}
\end{lemma}

\textbf{Proof:}

It suffices to consider cases where $H$ is positive (this is the only situation for which we need the lemma anyway). Since the oscillator functions are normalized:

\begin{equation*}
\left|\sqrt{n}\int_{-2}^{2} H(x,y) \psi_n(\sqrt{n}x)^2dx\right|\leq C
\end{equation*}

Fix $\epsilon\geq 0$. By the first two equalities of Lemma \ref{conv}, for $y$ fixed, we have:
\begin{equation*}
\lim_{n \to \infty}\int_{-2}^{2} \sqrt{n}H(x,y)\psi_n(\sqrt{n}x)^2dx=
\frac{1}{\pi}\int_{-2}^{2} H(x,y)\frac{1}{\sqrt{4-x^2}}dx
\end{equation*}

We can take this convergence to be uniform except on a set of measure $\frac{\epsilon}{8||f||_{\infty}}$, denoted $\Omega$. On $\Omega^c$, by the boundedness condition and the first two equalities of lemma \ref{conv} with an indicator function, we may select $N_1$ such that for $n\geq N_1$ we have:

\begin{equation*}
\left|\int^{2}_{-2}H(x,y)
\left(\sqrt{n}\psi_n(\sqrt{n}x)^2-\frac{1}{\pi\sqrt{4-x^2}}\right)dx\right|
\leq \frac{\epsilon}{4}
\end{equation*}

Therefore:
\begin{eqnarray*}
\Big|\int_{\Omega^c}\int_{-2}^{2} nH(x,y)\psi^2_n(\sqrt{n}x)\psi^2_n(\sqrt{n}y)dxdy\\
-\frac{1}{\pi}\int_{\Omega^c}\sqrt{n}\psi_n(\sqrt{n}y)^2
\int_{-2}^{2}H(x,y)\frac{1}{\sqrt{4-x^2}}dx\Big|\\
\leq \frac{\epsilon}{4}||\sqrt{n}\psi_n(\sqrt{n}y)^2||_1=\frac{\epsilon}{4}
\end{eqnarray*}

And, for $n\geq N_2$, for some sufficiently large $N_2$:
\begin{eqnarray*}
\Big|\int_{\Omega^c}\int_{-2}^{2} nH(x,y)\psi^2_n(\sqrt{n}x)\psi^2_n(\sqrt{n}y)dxdy\Big|\\
\leq
C\sqrt{n}\int_{-2}^21_{\Omega^c}(y)\psi^2_n(\sqrt{n}y)dy\leq C_f \frac{\epsilon}{4}
\end{eqnarray*}

In the last inequality, we have used the following:
\begin{eqnarray*}
\lim_{n \to \infty}\sqrt{n}\int_{-2}^21_{\Omega^c}(y)\psi^2_n(\sqrt{n}y)dy=
\int_{\Omega^c}\frac{1}{\pi\sqrt{4-x}}
\end{eqnarray*}

Consequently, for large enough $n$:
\begin{eqnarray*}
\Big|\int_{-2}^2\int_{-2}^{2} nH(x,y)\psi^2_n(\sqrt{n}x)\psi^2_n(\sqrt{n}y)dxdy\\
-\frac{1}{\pi}\int_{-2}^{2}\sqrt{n}\psi_n(\sqrt{n}y)^2
\int_{-2}^{2}H(x,y)\frac{1}{\sqrt{4-x^2}}dx\Big|
\leq \frac{\epsilon}{2}
\end{eqnarray*}

Applying the first two equalities of Lemma \ref{conv} again we obtain:
\begin{eqnarray*}
\Big|\frac{1}{\pi^2}\int_{-2}^{2}\int_{-2}^{2}\frac{1}{\sqrt{4-y^2}}
H(x,y)\frac{1}{\sqrt{4-x^2}}dx\\
-\frac{1}{\pi}\int_{-2}^{2}\sqrt{n}\psi_n(\sqrt{n}y)^2
\int_{-2}^{2}H(x,y)\frac{1}{\sqrt{4-x^2}}dx\Big|
\leq \frac{\epsilon}{2}
\end{eqnarray*}

The use of the first two equalities of Lemma \ref{conv} are justifed by $\int_{-2}^{2}H(x,y)\frac{1}{\sqrt{4-x^2}}dx$ being a bounded function of $y$, because one integrand is bounded and the second integrable. Consequently, for large enough $n$:

\begin{eqnarray*}
\Big|\int\int nH(x,y)\psi^2_n(\sqrt{n}x)\psi^2_n(\sqrt{n}y)dxdy-\\
\frac{1}{\pi^2}\int\int H(x,y)\frac{1}{\sqrt{4-x^2}\sqrt{4-y^2}}dxdy\Big|
\leq \epsilon
\end{eqnarray*}

Since $\epsilon$ was arbitrary, we have:
\begin{eqnarray*}
\lim_{n \to \infty}\int\int nH(x,y)\psi^2_n(\sqrt{n}x)\psi^2_n(\sqrt{n}y)dxdy
=\frac{1}{\pi^2}\int\int H(x,y)\frac{1}{\sqrt{4-x^2}\sqrt{4-y^2}}dxdy
\end{eqnarray*}

Arguing in the exact same way produces:
\begin{eqnarray*}
\lim_{n \to \infty}\int\int n H(x,y)\psi_n(\sqrt{n}x)\psi_{n-1}(\sqrt{n}x)\psi_n(\sqrt{n}y)
\psi_{n-1}(\sqrt{n}y)dxdy=\\
\frac{1}{4\pi^2}\int\int H(x,y)\frac{xy}{\sqrt{4-x^2}\sqrt{4-y^2}}dxdy
\end{eqnarray*}

Combining the two results concludes the proof. $\P$.

\subsection{Proof of Lemma \ref{boundbulk}}

First, we will aim to bound $\frac{1}{\sqrt{n}}\int_{-2}^{2} \psi^{\prime}_n(\sqrt{n}x)^{2}dx$.

Let $\mbox{Ai}(x)$ be the classical Airy function, which is defined by:
\begin{eqnarray*}
\mbox{Ai}(x)=\frac{1}{\pi}\lim_{b\to\infty}\int_0^b\cos\left(\frac{t^3}{3}+xt\right)dt
\end{eqnarray*}

By Lemma 3.7.2 in Zeitouni, we have the following identities on the very edge itself:
\begin{eqnarray*}
\lim_{n\to\infty}n^{1/12}\psi_n(2\sqrt{n})=\mbox{Ai}(0)\\
\lim_{n\to\infty}n^{-1/12}\psi^{'}_n(2\sqrt{n})=\mbox{Ai}^{\prime}(0)
\end{eqnarray*}
 Both of the values are finite, therefore $\psi_n(2\sqrt{n})\psi^{\prime}_n(2\sqrt{n})$ is convergent, and in particular is bounded.
Using integration by parts:

\begin{eqnarray*}
\frac{1}{\sqrt{n}}\int_{-2}^{2} \psi^{\prime}_n(\sqrt{n}x)^2dx
=-\frac{1}{\sqrt{n}}\int_{-2}^{2} \psi_n(\sqrt{n}x)\psi^{\prime \prime}_n(\sqrt{n}x)dx
+\frac{2}{n}\psi_n(2\sqrt{n})\psi^{\prime}_n(2\sqrt{n})\\
=\frac{1}{\sqrt{n}}\int_{-2}^{2} (n+\frac{1}{2}-\frac{nx^2}{4})\psi_n(\sqrt{n}x)^{2}dx
+\frac{2}{n}\psi_n(2\sqrt{n})\psi^{\prime}_n(2\sqrt{n})\\
\leq C\int_{-2}^{2} \sqrt{n}\psi_n(\sqrt{n}x)^{2}dx
+\frac{2}{n}\psi_n(2\sqrt{n})\psi^{\prime}_n(2\sqrt{n})
\end{eqnarray*}
This last equation is bounded since $\sqrt{n}\psi_n^{2}(\sqrt{n}x)$ has integral one, and we have proven that there exists a constant $C$ such that the following bound holds independently of $n$:
\begin{eqnarray*}
\frac{1}{\sqrt{n}}\int_{-2}^{2} \psi^{\prime}_n(\sqrt{n}x)^{2}dx\leq C
\end{eqnarray*}

We will next show that for any interval $I=[0,s]$ with $0 < s < 2$, there exists a constant $C_{s}$ such that $C_{s}\geq \sqrt{n}\psi_n(\sqrt{n}x)^2$ for all $x\in I$. We begin with the following equation, with for instance $C=2e^2$:

\begin{eqnarray*}
C\geq \sqrt{n}\int^{t}_{0}e^{x}\psi_n(\sqrt{n}x)^{2}dx=\sqrt{n}\int^{t}_{0}(e^{x})^{'}\psi_n(\sqrt{n}x)^{2}dx
\end{eqnarray*}

Notice that we have the trivial bound:
\begin{eqnarray*}
\sqrt{n}\int^{t}_{0}e^{x}\psi_n(\sqrt{n}x)^{2}dx \geq 0
\end{eqnarray*}

Integrating by parts:
\begin{eqnarray*}
C \geq -\sqrt{n}\int^{t}_{0}(e^{x})2\sqrt{n}\psi_n(\sqrt{n}x)\psi^{\prime}_n(\sqrt{n}x)dx
+e^{x}\sqrt{n}\psi^2_n(\sqrt{n}x)\Big|_{0}^t
\end{eqnarray*}

We next apply the following Hermite polynomial identity:
\begin{eqnarray*}
\psi_n(x)=\frac{\psi^{\prime \prime}_n(x)}{-n-\frac{1}{2}+\frac{x^2}{4}}
\end{eqnarray*}

This gives:
\begin{eqnarray*}
C \geq
\sqrt{n}\int^{t}_{0}\frac{e^{x}}{1/2+\frac{n}{4}(4-x^2)}2\sqrt{n}\psi_n(\sqrt{n}x)^{''}\psi^{\prime}_n(\sqrt{n}x)dx
+e^{x}\sqrt{n}\psi^2_n(\sqrt{n}x)\Big|_{0}^t
\end{eqnarray*}
Or else:
\begin{eqnarray*}
C \geq
\sqrt{n}\int^{t}_{0}\frac{e^{x}}{1/2+\frac{n}{4}(4-x^2)}(\psi_n^{\prime}(\sqrt{n}x)^2)^{'}dx
+e^{x}\sqrt{n}\psi^2_n(\sqrt{n}x)\Big|_{0}^t
\end{eqnarray*}
Integrating by parts again:
\begin{eqnarray*}
C +\sqrt{n}\int^{t}_{0}(\frac{e^{x}}{1/2+\frac{n}{4}(4-x^2)})^{'}(\psi_n^{'}(\sqrt{n}x)^2)dx\\
\geq \frac{\sqrt{n}e^{x}}{1/2+\frac{n}{4}(4-x^2)}\psi_n^{'}(\sqrt{n}t)^2\Big|_0^t
+e^{x}\sqrt{n}\psi^2_n(\sqrt{n}x)\Big|_{0}^t
\end{eqnarray*}

Expanding the derivative:
\begin{eqnarray*}
\sqrt{n}\left(\frac{e^{x}}{1/2+\frac{n}{4}(4-x^2)}\right)^{\prime}
=\sqrt{n}(\frac{e^{x}}{1/2+\frac{n}{4}(4-x^2)}+\frac{ne^{x}\frac{x}{2}}{(1/2+\frac{n}{4}(4-x^2))^2})\\
\leq
\sqrt{n}\left(
\frac{e^{x}}{\frac{n}{4}(4-x^2)}
+\frac{ne^{x}\frac{x}{2}}{\frac{n}{4}(4-x^2)}
\right)^2\\
\leq
\frac{1}{\sqrt{n}}\left(
\frac{e^{x}}{\frac{1}{4}(4-x^2)}
+\frac{e^{x}\frac{x}{2}}{\frac{1}{4}(4-x^2)}
\right)^2 \leq \frac{C}{\sqrt{n}}
\end{eqnarray*}

In the last inequality we have have used that $t$ is a bounded distance from the edge. Substituting this estimate into our previous expression, we have:
\begin{eqnarray*}
C \left(1+\frac{1}{\sqrt{n}}\int^{t}_{0}(\psi_n^{'}(\sqrt{n}x)^2)dx\right)\\
\geq \frac{\sqrt{n}e^{x}}{1/2+\frac{n}{4}(4-x^2)}\psi_n^{'}(\sqrt{n}x)^2\Big|_0^t
+e^{x}\sqrt{n}\psi^2_n(\sqrt{n}x)\Big|_{0}^t
\end{eqnarray*}

We have already seen that the right hand side is then bounded, and therefore we have:

\begin{eqnarray*}
C\geq \frac{\sqrt{n}e^{x}}{1/2+n/4(4-x^2)}\psi_n^{\prime}(\sqrt{n}x)^2 \Big|_0^t
+e^{x}\sqrt{n}\psi^2_n(\sqrt{n}x)\Big|_{0}^t
\end{eqnarray*}

And then, after discarding a positive term:
\begin{eqnarray*}
C+\frac{\sqrt{n}e^2}{1/2+n}\psi_n^{\prime}(0)^2+\sqrt{n}\psi^2_n(0) \\
\geq
\sqrt{n}\psi^2_n(\sqrt{n}t)
\end{eqnarray*}

By Stirling's formula and formulas (2.160) and (2.161) in \cite{Tao}, the right hand side is bounded by some constant:
\begin{eqnarray*}
C>\frac{\sqrt{n}e^2}{1/2+n}\psi_n^{\prime}(0)^2+\sqrt{n}\psi^2_n(0)
\end{eqnarray*}

This yields:

\begin{eqnarray*}
C\geq \sqrt{n}\psi^2_n(\sqrt{n}t)
\end{eqnarray*}

It remains to deal with the term $|\sqrt{n}\psi_n(\sqrt{n}t)\psi_{n-1}(\sqrt{n}t)|$, which we take care of presently. By a Hermite polynomial identity:

\begin{eqnarray*}
\psi_{n-1}(\sqrt{n}t)=\frac{1}{\sqrt{n}}\psi^{\prime}_n(\sqrt{n}t)+\frac{1}{2}\psi(\sqrt{n}t)
\end{eqnarray*}

By this identity, the Cauchy-Schwarz inequality, and our previous estimates we may compute:

\begin{eqnarray*}
|\sqrt{n}\psi_n(\sqrt{n}t)\psi_{n-1}(\sqrt{n}t)|
\leq \sqrt{n}\psi_n(\sqrt{n}t)^2+\sqrt{n}\psi_{n-1}(\sqrt{n}t)^2
\leq  C_s
\end{eqnarray*}

Consequently, for any $[0,s]\in(0,2)$ there exists a $C_s$ such that $|\sqrt{n}\psi_n(\sqrt{n}t)\psi_{n-1}(\sqrt{n}t)|\leq C_s$. We now have all the estimates which we sought, and the lemma is established. $\P$.

\end{document}